\documentclass[12pt,fleqn]{article}

\begin{document}

\title{Exponential Forms and Path Integrals for\\
Complex Numbers in $n$ Dimensions}

\author{Silviu Olariu
\thanks{e-mail: olariu@ifin.nipne.ro}\\
Institute of Physics and Nuclear Engineering, Tandem Laboratory\\
76900 Magurele, P.O. Box MG-6, Bucharest, Romania}

\date{28 July 2000}

\maketitle

\abstract

Two distinct systems of commutative complex numbers in $n$ dimensions are
described, of polar and planar types.  Exponential forms of n-complex
numbers are given in each case, which depend on geometric variables. Azimuthal 
angles, which are cyclic variables, appear in these forms at the exponent, 
and this leads to the concept of residue for path integrals of n-complex
functions. The exponential function of an n-complex number is expanded in
terms of functions called in this paper cosexponential functions,
which are generalizations to $n$ dimensions of the circular and hyperbolic sine
and cosine functions.  The factorization of n-complex polynomials is discussed.

\endabstract

\section{Introduction}

Hypercomplex numbers are a generalization to several dimensions
of the regular complex numbers in 2 dimensions. 
A well-known example
of hypercomplex numbers are the quaternions of Hamilton, which are a system of
hypercomplex numbers in 
four dimensions, the multiplication being a non-commutative operation. \cite{1}
Many other hypercomplex systems are possible, \cite{2a}-\cite{2b}
but these systems do not have all the required properties of regular,
two-dimensional complex numbers which rendered possible the development of the
theory of functions of a 2-dimensional complex variable.

Two distinct systems of complex numbers in $n$ dimensions are
described in this paper, for which the multiplication is associative and
commutative, and which are rich enough in properties such that exponential
forms 
exist and the concepts of analytic n-complex function,  contour integration and
residue can be defined. 
The first type of n-complex numbers described in this article is
characterized by the presence, 
in an odd number of dimensions, of one polar axis, and by the presence, 
in an even number of dimensions, of two polar axes.
Therefore, these numbers will be called polar n-complex numbers. 
The other type of n-complex numbers described in this paper exists as a
distinct entity only when the number of dimensions $n$ of the space is even.
These numbers will be called planar n-complex numbers. 
The planar hypercomplex numbers become for $n=2$ the
usual complex numbers $x+iy$.

The central result of this paper is the existence of an exponential form of
n-complex numbers, which is expressed in terms of geometric variables. 
The
exponential form provides the link between the algebraic side of
the operations and the analytic properties of the functions of n-complex
variables. The azimuthal angles 
$\phi_k$, which are cyclic variables, appear in these forms at the exponent,
and this leads to the concept of  n-complex residue for path integrals
of n-complex functions. Expressions are 
given for the elementary functions of n-complex variable.  The
exponential function of an n-complex number is expanded in terms of functions
called in this paper n-dimensional cosexponential functions of the polar and
respectively planar type.  The polar cosexponential functions are a
generalization to $n$ dimensions of the hyperbolic functions $\cosh y, \sinh
y$, and the planar cosexponential functions are a generalization to $n$
dimensions of the trigonometric functions $\cos y, \sin y$. Addition theorems
and other relations are obtained for the n-dimensional cosexponential
functions.

In the case of polar n-complex numbers, a polynomial 
can be written as a product of linear or
quadratic factors, although several factorizations are in general possible.  In
the case of planar n-complex numbers, a polynomial 
can always be written as a product of
linear factors, although, again, several factorizations are in general
possible.

\section{Polar n-complex numbers}

\subsection{Operations with polar n-complex numbers}

A hypercomplex number in $n$ dimensions is determined by its $n$ components
$(x_0,x_1,...,x_{n-1})$. The polar n-complex numbers and
their operations discussed in this paper can be represented 
by  writing the n-complex number $(x_0,x_1,...,x_{n-1})$ as  
$u=x_0+h_1x_1+h_2x_2+\cdots+h_{n-1}x_{n-1}$, where 
$h_1, h_2, \cdots, h_{n-1}$ are bases for which the multiplication rules are 
\begin{equation}
h_j h_k =h_l ,\:l=j+k-n[(j+k)/n],
\label{1}
\end{equation}
for $j,k,l=0,1,..., n-1$, where $h_0=1$.
In this relation, $[(j+k)/n]$ denotes the integer part of $(j+k)/n$, defined as
$[a]\leq a<[a]+1$, so that $0\leq j+k-n[(j+k)/n]\leq n-1$. 
In this paper, brackets larger than the regular brackets
$[\;]$ do not have the meaning of integer part.
The significance of the composition laws in Eq.
(\ref{1}) can be understood by representing the bases $h_j, h_k$ by points on a
circle at the angles $\alpha_j=2\pi j/n,\alpha_k=2\pi k/n$, as shown in Fig. 1,
and the product $h_j h_k$ by the point of the circle at the angle 
$2\pi (j+k)/n$. If $2\pi\leq 2\pi (j+k)/n<4\pi$, the point represents the basis
$h_l$ of angle $\alpha_l=2\pi(j+k)/n-2\pi$.

Two n-complex numbers 
$u=x_0+h_1x_1+h_2x_2+\cdots+h_{n-1}x_{n-1}$,
$u^\prime=x^\prime_0+h_1x^\prime_1+h_2x^\prime_2+\cdots+h_{n-1}x^\prime_{n-1}$ 
are equal if and only if $x_j=x^\prime_j, j=0,1,...,n-1$.
The sum of the n-complex numbers $u$
and
$u^\prime$ 
is
\begin{equation}
u+u^\prime=x_0+x^\prime_0+h_1(x_1+x^\prime_1)+\cdots
+h_{n-1}(x_{n-1} +x^\prime_{n-1}) .
\label{2}
\end{equation}
The product of the n-complex numbers $u, u^\prime$ is 
\begin{equation}
\begin{array}{l}
uu^\prime=x_0 x_0^\prime +x_1x_{n-1}^\prime+x_2 x_{n-2}^\prime
+x_3x_{n-3}^\prime+\cdots+x_{n-1}x_1^\prime\\
+h_1(x_0 x_1^\prime+x_1x_0^\prime+x_2x_{n-1}^\prime+x_3x_{n-2}^\prime
+\cdots+x_{n-1} x_2^\prime) \\
+h_2(x_0 x_2^\prime+x_1x_1^\prime+x_2x_0^\prime+x_3x_{n-1}^\prime
+\cdots+x_{n-1} x_3^\prime) \\
\vdots\\
+h_{n-1}(x_0 x_{n-1}^\prime+x_1x_{n-2}^\prime+x_2x_{n-3}^\prime
+x_3x_{n-4}^\prime
+\cdots+x_{n-1} x_0^\prime).
\end{array}
\label{3}
\end{equation}
The product $uu^\prime$ can be written as
\begin{equation}
uu^\prime=\sum_{k=0}^{n-1}h_k\sum_{l=0}^{n-1}x_l x^\prime_{k-l+n[(n-k-1+l)/n]}.
\label{3a}
\end{equation}
If $u,u^\prime,u^{\prime\prime}$ are n-complex numbers, the multiplication 
is associative
$(uu^\prime)u^{\prime\prime}=u(u^\prime u^{\prime\prime})$
and commutative
$u u^\prime=u^\prime u ,$
because the product of the bases, defined in Eq. (\ref{1}), 
is associative and commutative. 

The inverse of polar the n-complex number $u$ is the n-complex number
$u^\prime$ having the property that $uu^\prime=1$.
This equation has a solution provided that the corresponding determinant
$\nu$ is not equal to zero, $\nu\not=0$.
If $n$ is an even number, it can be shown that
\begin{equation}
\nu=v_+v_-\prod_{k=1}^{n/2-1}\rho_k^2,
\end{equation}
and if $n$ is an odd number,
\begin{equation}
\nu=v_+\prod_{k=1}^{(n-1)/2}\rho_k^2, 
\end{equation}
where 
\begin{equation}
\rho_k^2= v_k^2+\tilde v_k^2,   
\end{equation}
\begin{equation}
v_k=\sum_{p=0}^{n-1}x_p\cos\left(\frac{2\pi kp}{n}\right),
\end{equation}
\begin{equation}
\tilde v_k=\sum_{p=0}^{n-1}x_p\sin\left(\frac{2\pi kp}{n}\right).
\end{equation}
Thus, in an even number of dimensions $n$, an n-complex number has an inverse
unless it lies on one of the nodal hypersurfaces $v_+=0$, or
$v_-=0$, or $\rho_1=0$, or ... or $\rho_{n/2-1}=0$.
In an odd number of dimensions $n$, an n-complex number has an inverse
unless it lies on one of the nodal hypersurfaces $v+=0$,
or $\rho_1=0$, or ... or $\rho_{(n-1)/2}=0$. 

For even $n$,
\begin{equation}
d^2=\frac{1}{n}v_+^2+\frac{1}{n}v_-^2
+\frac{2}{n}\sum_{k=1}^{n/2-1}\rho_k^2,
\end{equation}
and for odd $n$,
\begin{equation}
d^2=\frac{1}{n}v_+^2
+\frac{2}{n}\sum_{k=1}^{(n-1)/2}\rho_k^2 .
\end{equation}
From these relations it results that if the product of two n-complex numbers is
zero, $uu^\prime=0$, then 
$\rho_+\rho_+^\prime=0, \rho_-\rho_-^\prime=0,
\rho_k\rho_k^\prime=0, k=1,...,n/2$, which means that either $u=0$, or
$u^\prime=0$, or $u, u^\prime$ belong to orthogonal hypersurfaces in such a way
that the afore-mentioned products of components should be equal to zero.

\subsection{Geometric representation of polar n-complex numbers}

The polar n-complex number $x_0+h_1x_1+h_2x_2+\cdots+h_{n-1}x_{n-1}$
can be represented by 
the point $A$ of coordinates $(x_0,x_1,...,x_{n-1})$. 
If $O$ is the origin of the n-dimensional space,  the
distance from the origin $O$ to the point $A$ of coordinates
$(x_0,x_1,...,x_{n-1})$ has the expression
\begin{equation}
d^2=x_0^2+x_1^2+\cdots+x_{n-1}^2.
\label{10}
\end{equation}
The quantity $d$ will be called modulus of the polar n-complex number 
$u=x_0+h_1x_1+h_2x_2+\cdots+h_{n-1}x_{n-1}$. The modulus of an n-complex number
$u$ will be designated by $d=|u|$. If $\nu>0$, the quantity $\rho=\nu^{1/n}$
will be called amplitude of the n-complex number $u$. 

The exponential and trigonometric forms of the n-complex number $u$ can be
obtained conveniently in a rotated system of axes defined by the transformation
\begin{equation}
v_+= \sqrt{n}\xi_+ ,  v_-=  \sqrt{n}\xi_-,  
v_k= \sqrt{n/2}\xi_k , \tilde v_k= \sqrt{n/2}\eta_k, 
\end{equation}
for $k=1,...,[(n-1)/2]$ .
This transformation from the coordinates $x_0,...,x_{n-1}$ to the variables
$\xi_+,\xi_-, \xi_k, \eta_k$ is unitary.

The position of the point $A$ of coordinates $(x_0,x_1,...,x_{n-1})$ can be
also described with the aid of the distance $d$, Eq. (\ref{10}), and of $n-1$
angles defined further. Thus, in the plane of the axes $v_k,\tilde v_k$, the
azimuthal angles $\phi_k$ can be introduced by the relations 
\begin{equation}
\cos\phi_k=v_k/\rho_k,\:\sin\phi_k=\tilde v_k/\rho_k, 
\end{equation}
where $0\leq\phi_k<2\pi$, 
so that there are $[(n-1)/2]$ azimuthal angles.
If the projection of the point $A$ on the plane of the axes $v_k,\tilde v_k$ is
$A_k$, 
and the projection of the point $A$ on the 4-dimensional space defined by the
axes $v_1, \tilde v_1, v_k,\tilde v_k$ is $A_{1k}$, the angle $\psi_{k-1}$
between the line 
$OA_{1k}$ and the 2-dimensional plane defined by the axes $v_k,\tilde v_k$ is 
\begin{equation}
\tan\psi_{k-1}=\rho_1/\rho_k, 
\end{equation}
for $0\leq\psi_k\leq\pi/2, k=2,...,[(n-1)/2]$,
so that there are $[(n-3)/2]$ planar angles.
Moreover, there is a polar angle $\theta_+$, which can be
defined as the angle between the line $OA_{1+}$ and the axis $v_+$,
where $A_{1+}$ is the projection of the point $A$ on the 3-dimensional space
generated by the axes $v_1, \tilde v_1, v_+$,
\begin{equation}
\tan\theta_+=\sqrt{2}\rho_1/v_+, 
\end{equation}
where $0\leq\theta_+\leq\pi$ , 
and in an even number of dimensions $n$ there is also a polar angle $\theta_-$,
which can be defined as the angle between the line $OA_{1-}$ and the axis
$v_-$, 
where $A_{1-}$ is the projection of the point $A$ on the 3-dimensional space
generated by the axes $v_1, \tilde v_1, v_-$,
\begin{equation}
\tan\theta_-=\sqrt{2}\rho_1/v_-, 
\end{equation}
where $0\leq\theta_-\leq\pi$ .
Thus, the position of
the point $A$ is described, in an even number of dimensions, by the distance
$d$, by $n/2-1$ azimuthal angles, by $n/2-2$ planar angles, and by 2 polar
angles. In an odd number of dimensions, the position of the point $A$ is
described by $(n-1)/2$ azimuthal angles, by $(n-3)/2$ planar angles, and by 1
polar angle. These angles are shown in Fig.  2. The variables $\nu, \rho,
\rho_k, \tan\theta_+/\sqrt{2}, \tan\theta_-/\sqrt{2}, 
\tan\psi_k$ are multiplicative and the azimuthal angles $\phi_k$ are additive
upon the multiplication of polar n-complex numbers.

\subsection{The n-dimensional polar cosexponential functions}

The exponential function of the polar n-complex variable $u$ can be defined by
the series
$\exp u = 1+u+u^2/2!+u^3/3!+\cdots$ . 
It can be checked by direct multiplication of the series that
$\exp(u+u^\prime)=\exp u \cdot \exp u^\prime$ ,
so that  
$\exp u=\exp x_0 \cdot \exp (h_1x_1) \cdots \exp (h_{n-1}x_{n-1})$.

It can be seen with the aid of the representation in Fig. 1 that 
\begin{equation}
h_k^{n+p}=h_k^p, 
\end{equation}
for $p$ integer,  $k=1,...,n-1$.
Then $e^{h_k y}$ can be written as
\begin{equation}
e^{h_k y}=\sum_{p=0}^{n-1}h_{kp-n[kp/n]}g_{nl}(y),
\label{28b}
\end{equation}
where the functions $g_{nl}$, which will be
called polar cosexponential functions in $n$ dimensions, are
\begin{equation}
g_{nl}(y)=\sum_{p=0}^\infty \frac{y^{l+pn}}{(l+pn)!},
\label{29}
\end{equation}
for $l=0,1,...,n-1$.
If $n$ is even, the polar cosexponential functions of even index $k$ are
even functions, $g_{n,2l}(-y)=g_{n,2l}(y)$, 
and the polar cosexponential functions of odd index 
are odd functions, $g_{n,2l+1}(-y)=-g_{n,2l+1}(y)$, $l=0,1,...,n/2-1$. For odd
values of $n$, the polar cosexponential functions do not have a definite
parity. It can be checked that
\begin{equation}
\sum_{l=0}^{n-1}g_{nl}(y)=e^y
\end{equation}
and, for even $n$, 
\begin{equation}
\sum_{l=0}^{n-1}(-1)^k g_{nl}(y)=e^{-y}.
\end{equation}
The expression of the polar n-dimensional cosexponential functions is
\begin{equation}
g_{nk}(y)
=\frac{1}{n}\sum_{l=0}^{n-1}\exp\left[y\cos\left(\frac{2\pi l}{n}\right)\right]
\cos\left[y\sin\left(\frac{2\pi l}{n}\right)-\frac{2\pi kl}{n}\right],
\label{30}
\end{equation}
for $k=0,1,...,n-1$.
It can be shown from Eq. (\ref{30}) that
\begin{equation}
\sum_{k=0}^{n-1}g_{nk}^2(y)=\frac{1}{n}\sum_{l=0}^{n-1}\exp\left[2y\cos\left(
\frac{2\pi l}{n}\right)\right].
\label{34a}
\end{equation}
It can be seen that the right-hand side of Eq. (\ref{34a}) does not contain
oscillatory terms. 
If $n$ is a multiple of 4, it can be shown by replacing $y$
by $iy$ in Eq. (\ref{34a}) that
\begin{equation}
\sum_{k=0}^{n-1}(-1)^kg_{nk}^2(y)=\frac{2}{n}\left\{1+\cos 2y
+\sum_{l=1}^{n/4-1}\cos\left[2y\cos\left(\frac{2\pi l}{n}\right)\right]
\right\},
\label{34b}
\end{equation}
which does not contain exponential terms.

Addition theorems for the polar n-dimensional cosexponential functions can be
obtained from the relation $\exp h_1(y+z)=\exp h_1 y \cdot\exp h_1 z $, by
substituting the expression of the exponentials as given 
by $e^{h_1 y}=\sum_{p=0}^{n-1} h_p g_{np}(y)$,
\begin{eqnarray}
\lefteqn{g_{nk}(y+z)=g_{n0}(y)g_{nk}(z)+g_{n1}(y)g_{n,k-1}(z)
+\cdots+g_{nk}(y)g_{n0}(z)\nonumber}\\
&&+g_{n,k+1}(y)g_{n,n-1}(z)+g_{n,k+2}(y)g_{n,n-2}(z)+\cdots
+g_{n,n-1}(y)g_{n,k+1}(z) ,
\label{35a}
\end{eqnarray}
for $k=0,1,...,n-1$.

It can also be shown that
\begin{equation}
\left\{g_{n0}(y)+h_1g_{n1}(y)+\cdots+h_{n-1}g_{n,n-1}(y)\right\}^l
=g_{n0}(ly)+h_1g_{n1}(ly)+\cdots+h_{n-1}g_{n,n-1}(ly).
\label{37b}
\end{equation}

The polar n-dimensional cosexponential functions are solutions of the
$n^{\rm th}$-order differential equation
\begin{equation}
\frac{{\rm d}^n\zeta}{{\rm d}u^n}=\zeta ,
\label{44}
\end{equation}
whose solutions are of the form
$\zeta(u)=A_0g_{n0}(u)+A_1g_{n1}(u)+\cdots+A_{n-1}g_{n,n-1}(u)$.
It can be checked that the derivatives of the polar cosexponential functions
are related by
\begin{equation}
\frac{dg_{n0}}{du}=g_{n,n-1}, \:
\frac{dg_{n1}}{du}=g_{n0}, \:...,
\frac{dg_{n,n-2}}{du}=g_{n,n-3} ,
\frac{dg_{n,n-1}}{du}=g_{n,n-2} .
\label{45}
\end{equation}

For $n=2$, the polar cosexponential functions are $g_{20}(y)=\cosh y$ and
$g_{21}(y)=\sinh y$.

\subsection{Exponential and trigonometric forms of polar n-complex numbers}

In order to obtain the exponential and trigonometric forms of polar n-complex
numbers, a new set of hypercomplex bases 
will be introduced for even $n$
by the relations 
\begin{equation}
e_+=\frac{1}{n}\sum_{p=0}^{n-1}h_p, 
\label{e13a}
\end{equation}
\begin{equation}
e_k=\frac{2}{n}\sum_{p=0}^{n-1}h_p\cos\left(\frac{2\pi kp}{n}\right),
\end{equation}
\begin{equation}
\tilde e_k=\frac{2}{n}\sum_{p=0}^{n-1}h_p\sin\left(\frac{2\pi kp}{n}\right),
\end{equation}
where $k=1,...,[(n-1)/2]$ and, if $n$ is even,  
\begin{equation}
e_-=\frac{1}{n}\sum_{p=0}^{n-1}(-1)^p h_p.
\label{e13b}
\end{equation}
The multiplication relations for the new hypercomplex bases are
\begin{eqnarray}
\lefteqn{e_+^2=e_+,\; e_-^2=e_-,\; e_+e_-=0,\; e_+e_k=0,\; e_+\tilde e_k=0,\;
e_-e_k=0,\; 
e_-\tilde e_k=0,\nonumber}\\ 
&&e_k^2=e_k,\; \tilde e_k^2=-e_k,\; e_k \tilde e_k=\tilde e_k ,\; e_ke_l=0,\;
e_k\tilde e_l=0,\; \tilde e_k\tilde e_l=0,\; k\not=l,\; 
\label{e12a}
\end{eqnarray}
where $k,l=1,...,[(n-1)/2]$.
It can be shown that, for even $n$,
\begin{equation}
u=e_+v_+ + e_-v_- 
+\sum_{k=1}^{n/2-1}(e_k v_k+\tilde e_k \tilde v_k),
\end{equation}
and for odd $n$
\begin{equation}
u=e_+v_+ +\sum_{k=1}^{(n-1)/2}(e_k v_k+\tilde e_k \tilde v_k).
\end{equation}

The exponential form of the n-complex number $u$ is
\begin{eqnarray}\lefteqn{
u=\rho\exp\left\{\sum_{p=1}^{n-1}h_p\left[
\frac{1}{n}\ln\frac{\sqrt{2}}{\tan\theta_+}
+F(n)\frac{(-1)^p}{n}\ln\frac{\sqrt{2}}{\tan\theta_-}\right.\right.\nonumber}\\
&&
\left.\left.-\frac{2}{n}\sum_{k=2}^{[(n-1)/2]}
\cos\left(\frac{2\pi kp}{n}\right)\ln\tan\psi_{k-1}
\right]
+\sum_{k=1}^{[(n-1)/2]}\tilde e_k\phi_k ,
\right\},
\label{50a}
\end{eqnarray}
where $F(n)=1$ for even $n$ and $F(n)=0$ for odd $n$, and 
\begin{equation}
\rho=\left(v_+v_-\rho_1^2\cdots \rho_{n/2-1}^2\right)^{1/n}
\end{equation}
for even $n$, and
\begin{equation}
\rho=\left(v_+\rho_1^2\cdots \rho_{(n-1)/2}^2\right)^{1/n}
\end{equation}
for odd $n$.

The trigonometric form of the n-complex number $u$ is
\begin{eqnarray}
\lefteqn{u=d
\left(\frac{n}{2}\right)^{1/2}
\left(\frac{1}{\tan^2\psi_+}+\frac{F(n)}{\tan^2\psi_-}+1
+\frac{1}{\tan^2\psi_1}+\frac{1}{\tan^2\psi_2}+\cdots
+\frac{1}{\tan^2\psi_{[(n-3)/2]}}\right)^{-1/2}\nonumber}\\
&&\left(\frac{e_+\sqrt{2}}{\tan\theta_+}+F(n)\frac{e_-\sqrt{2}}{\tan\theta_-}
+e_1+\sum_{k=2}^{[(n-1)/2]}\frac{e_k}{\tan\psi_{k-1}}\right)
\exp\left(\sum_{k=1}^{[(n-1)/2]}\tilde e_k\phi_k\right).
\label{52a}
\end{eqnarray}

\subsection{Elementary functions of a polar n-complex variable}

The logarithm $u_1$ of the polar n-complex number $u$, $u_1=\ln u$, can be
defined 
as the solution of the equation
$u=e^{u_1}$ .
For even $n$, $\ln u$ exists 
as an n-complex function with real components if $v_+>0$
and $v_->0$.
For odd $n$ $\ln u$ exists 
as an n-complex function with real components if
$v_+>0$. 
The expression of the logarithm is
\begin{eqnarray}\lefteqn{
\ln u=\ln \rho+\sum_{p=1}^{n-1}h_p\left[
\frac{1}{n}\ln\frac{\sqrt{2}}{\tan\theta_+}
+F(n)\frac{(-1)^p}{n}\ln\frac{\sqrt{2}}{\tan\theta_-}\right.
\left.-\frac{2}{n}\sum_{k=2}^{[(n-1)/2]}
\cos\left(\frac{2\pi kp}{n}\right)\ln\tan\psi_{k-1}
\right]\nonumber}\\
&&+\sum_{k=1}^{[(n-1)/2]}\tilde e_k\phi_k.
\label{56a}
\end{eqnarray}
The function $\ln u$ is multivalued because of the presence of the terms 
$\tilde e_k\phi_k$.

The power function $u^m$ of the polar n-complex variable $u$ can be defined for
real values of $m$ as $u^m=e^{m\ln u}$ .
It can be shown that
\begin{equation}
u^m=e_+ v_+^m+F(n)e_- v_-^m +\sum_{k=1}^{[(n-1)/2]}
\rho_k^m(e_k\cos m\phi_k+\tilde e_k\sin m\phi_k).
\label{59a}
\end{equation}
For integer values of $m$, this expression is valid 
for any $x_0,...,x_{n-1}$.
The power function is multivalued unless $m$ is an integer.

\subsection{Power series of polar n-complex numbers}

A power series of the polar n-complex variable $u$ is a series of the form
\begin{equation}
a_0+a_1 u + a_2 u^2+\cdots +a_l u^l+\cdots .
\label{83}
\end{equation}
Using the inequality 
\begin{equation}
|u^\prime u^{\prime\prime}|\leq
\sqrt{n}|u^\prime||u^{\prime\prime}| , 
\label{83b}
\end{equation}
which replaces the relation of equality
extant for 2-dimensional complex numbers, it can be shown that
the series (\ref{83}) is absolutely convergent for 
$|u|<c$,
where 
$c=\lim_{l\rightarrow\infty} |a_l|/\sqrt{n}|a_{l+1}|$ .

The convergence of the series (\ref{83}) can be also studied with the aid of
the formulas (\ref{59a}), which for integer values of $m$ are
valid for any values of $x_0,...,x_{n-1}$.
If $a_l=\sum_{p=0}^{n-1}h_p a_{lp}$, and
\begin{equation}
A_{l+}=\sum_{p=0}^{n-1}a_{lp},
\label{88a}
\end{equation}
\begin{equation}
A_{lk}=\sum_{p=0}^{n-1}a_{lp}\cos\left(\frac{2\pi kp}{n}\right),
\label{88b}
\end{equation}
\begin{equation}
\tilde A_{lk}=\sum_{p=0}^{n-1}a_{lp}\sin\left(\frac{2\pi kp}{n}\right),
\label{88c}
\end{equation}
for $k=1,...,[(n-1)/2]$, and for even $n$ 
\begin{equation}
A_{l-}=\sum_{p=0}^{n-1}(-1)^p a_{lp},
\label{88d}
\end{equation}
the series (\ref{83}) can be written as
\begin{equation}
\sum_{l=0}^\infty \left[
e_+A_{l+}v_+^l+F(n)e_-A_{l-}v_-^l+\sum_{k=1}^{[(n-1)/2]}
(e_k A_{lk}+\tilde e_k\tilde A_{lk})(e_k v_k+\tilde e_k\tilde v_k)^l 
\right].
\label{89a}
\end{equation}

The series in Eq. (\ref{83}) is absolutely convergent for 
\begin{equation}
|v_+|<c_+,\:
|v_-|<c_-,\:
\rho_k<c_k, 
\end{equation}
for $k=1,..., [(n-1)/2]$,
where
\begin{equation}
c_+=\lim_{l\rightarrow\infty} \frac{|A_{l+}|}{|A_{l+1,+}|} ,\:
\end{equation}
\begin{equation}
c_-=\lim_{l\rightarrow\infty} \frac{|A_{l-}|}{|A_{l+1,-}|} ,\:
\end{equation}
\begin{equation}
c_k=\lim_{l\rightarrow\infty} \frac
{\left(A_{lk}^2+\tilde A_{lk}^2\right)^{1/2}}
{\left(A_{l+1,k}^2+\tilde A_{l+1,k}^2\right)^{1/2}} .
\end{equation}
These relations show that the region of convergence of the series
(\ref{83}) is an n-dimensional cylinder.

\subsection{Analytic functions of polar n-complex variables}

The derivative  
of a function $f(u)$ of the n-complex variables $u$ is
defined as a function $f^\prime (u)$ having the property that
\begin{equation}
|f(u)-f(u_0)-f^\prime (u_0)(u-u_0)|\rightarrow 0 \:\:{\rm as} 
\:\:|u-u_0|\rightarrow 0 . 
\label{h88}
\end{equation}
If the difference $u-u_0$ is not parallel to one of the nodal hypersurfaces,
the definition in Eq. (\ref{h88}) can also 
be written as
\begin{equation}
f^\prime (u_0)=\lim_{u\rightarrow u_0}\frac{f(u)-f(u_0)}{u-u_0} .
\label{h89}
\end{equation}
The derivative of the function $f(u)=u^m $, with $m$ an integer, 
is $f^\prime (u)=mu^{m-1}$, as can be seen by developing $u^m=[u_0+(u-u_0)]^m$
as
\begin{equation}
u^m=\sum_{p=0}^{m}\frac{m!}{p!(m-p)!}u_0^{m-p}(u-u_0)^p,
\label{h90}
\end{equation}
and using the definition (\ref{h88}).

If the function $f^\prime (u)$ defined in Eq. (\ref{h88}) is independent of the
direction in space along which $u$ is approaching $u_0$, the function $f(u)$ 
is said to be analytic, analogously to the case of functions of regular complex
variables. \cite{3} 
The function $u^m$, with $m$ an integer, 
of the n-complex variable $u$ is analytic, because the
difference $u^m-u_0^m$ is always proportional to $u-u_0$, as can be seen from
Eq. (\ref{h90}). Then series of
integer powers of $u$ will also be analytic functions of the n-complex
variable $u$, and this result holds in fact for any commutative algebra. 

If the n-complex function $f(u)$
of the polar n-complex variable $u$ is written in terms of 
the real functions $P_k(x_0,...,x_{n-1}), k=0,1,...,n-1$ of the real
variables $x_0,x_1,...,x_{n-1}$ as 
\begin{equation}
f(u)=\sum_{k=0}^{n-1}h_kP_k(x_0,...,x_{n-1}),
\label{h93}
\end{equation}
then relations of equality 
exist between the partial derivatives of the functions $P_k$. 
The derivative of the function $f$ can be written as
\begin{eqnarray}
\lim_{\Delta u\rightarrow 0}\frac{1}{\Delta u} 
\sum_{k=0}^{n-1}\left(h_k\sum_{l=0}^{n-1}
\frac{\partial P_k}{\partial x_l}\Delta x_l\right),
\label{h94}
\end{eqnarray}
where
$\Delta u=\sum_{k=0}^{n-1}h_l\Delta x_l$.

The
relations between the partials derivatives of the functions $P_k$ are
obtained by setting successively in   
Eq. (\ref{h94}) $\Delta u=h_l\Delta x_l$, for $l=0,1,...,n-1$, and equating the
resulting expressions. 
The relations are 
\begin{equation}
\frac{\partial P_k}{\partial x_0} = \frac{\partial P_{k+1}}{\partial x_1} 
=\cdots=\frac{\partial P_{n-1}}{\partial x_{n-k-1}} 
= \frac{\partial P_0}{\partial x_{n-k}}=\cdots
=\frac{\partial P_{k-1}}{\partial x_{n-1}}, 
\label{h95}
\end{equation}
for $ k=0,1,...,n-1$.
The relations (\ref{h95}) are analogous to the Riemann relations
for the real and imaginary components of a complex function. 
It can be shown from Eqs. (\ref{h95}) that the components $P_k$ fulfil the
second-order equations
\begin{eqnarray}
\lefteqn{\frac{\partial^2 P_k}{\partial x_0\partial x_l}
=\frac{\partial^2 P_k}{\partial x_1\partial x_{l-1}}
=\cdots=
\frac{\partial^2 P_k}{\partial x_{[l/2]}\partial x_{l-[l/2]}}}\nonumber\\
&&=\frac{\partial^2 P_k}{\partial x_{l+1}\partial x_{n-1}}
=\frac{\partial^2 P_k}{\partial x_{l+2}\partial x_{n-2}}
=\cdots
=\frac{\partial^2 P_k}{\partial x_{l+1+[(n-l-2)/2]}
\partial x_{n-1-[(n-l-2)/2]}} ,
\label{96}
\end{eqnarray}
for $k,l=0,1,...,n-1$.

\subsection{Integrals of polar n-complex functions}

The singularities of polar n-complex functions arise from terms of the form
$1/(u-u_0)^m$, with $m>0$. Functions containing such terms are singular not
only at $u=u_0$, but also at all points of the hypersurfaces
passing through $u_0$ and which are parallel to the nodal hypersurfaces. 

The integral of a polar n-complex function between two points $A, B$ along a
path 
situated in a region free of singularities is independent of path, which means
that the integral of an analytic function along a loop situated in a region
free of singularities is zero,
\begin{equation}
\oint_\Gamma f(u) du = 0,
\label{111}
\end{equation}
where it is supposed that a surface $\Sigma$ spanning 
the closed loop $\Gamma$ is not intersected by any of
the hypersurfaces associated with the
singularities of the function $f(u)$. Using the expression, Eq. (\ref{h93}),
for $f(u)$ and the fact that 
$du=\sum_{k=0}^{n-1}h_k dx_k$, 
the explicit form of the integral in Eq. (\ref{111}) is
\begin{equation}
\oint _\Gamma f(u) du = \oint_\Gamma
\sum_{k=0}^{n-1}h_k\sum_{l=0}^{n-1}P_l dx_{k-l+n[(n-k-1+l)/n]}.
\end{equation}

If the functions $P_k$ are regular on a surface $\Sigma$
spanning the loop $\Gamma$,
the integral along the loop $\Gamma$ can be transformed in an integral over the
surface $\Sigma$ of terms of the form
$\partial P_l/\partial x_{k-m+n[(n-k+m-1)/n]} $
$-\partial P_m/\partial x_{k-l+n[(n-k+l-1)/n]}$.
The integrals of these terms are equal to zero by Eqs. (\ref{h95}), and this
proves Eq. (\ref{111}). 

The quantity $du/(u-u_0)$ is
\begin{eqnarray}\lefteqn{
\frac{du}{u-u_0}=
\frac{d\rho}{\rho}
+\sum_{p=1}^{n-1}h_p\left[
\frac{1}{n}d\ln\frac{\sqrt{2}}{\tan\theta_+}
+F(n)\frac{(-1)^p}{n}d\ln\frac{\sqrt{2}}{\tan\theta_-}\right.\nonumber}\\
&&
\left.-\frac{2}{n}\sum_{k=2}^{[(n-1)/2]}
\cos\left(\frac{2\pi kp}{n}\right)d\ln\tan\psi_{k-1}
\right]
+\sum_{k=1}^{[(n-1)/2]}\tilde e_kd\phi_k.
\label{114a}
\end{eqnarray}
Since $\rho, \ln(\sqrt{2}/\tan\theta_+),\ln(\sqrt{2}/\tan\theta_-),
\ln(\tan\psi_{k-1})$ are singlevalued variables, it follows that
$\oint_\Gamma d\rho/\rho =0, 
\oint_\Gamma d(\ln\sqrt{2}/\tan\theta_+)=0,
\oint_\Gamma d(\ln\sqrt{2}/\tan\theta_-)=0,
\oint_\Gamma d(\ln\tan\psi_{k-1})=0$.
On the other hand, since
$\phi_k$ are cyclic variables, they may give contributions to
the integral around the closed loop $\Gamma$.

The expression of $\oint_\Gamma du/(u-u_0)$ can be written 
with the aid of a functional which will be called int($M,C$), defined for a
point $M$ and a closed curve $C$ in a two-dimensional plane, such that 
\begin{equation}
{\rm int}(M,C)=\left\{
\begin{array}{l}
1 \;\:{\rm if} \;\:M \;\:{\rm is \;\:an \;\:interior \;\:point \;\:of} \;\:C
,\\  
0 \;\:{\rm if} \;\:M \;\:{\rm is \;\:exterior \;\:to}\:\; C .\\
\end{array}\right.
\label{118}
\end{equation}
If $f(u)$ is an analytic function of a polar n-complex variable which can be
expanded in a 
series in the region of the curve
$\Gamma$ and on a surface spanning $\Gamma$, then 
\begin{equation}
\oint_\Gamma \frac{f(u)du}{u-u_0}=
2\pi f(u_0)\sum_{k=1}^{[(n-1)/2]}\tilde e_k 
\;{\rm int}(u_{0\xi_k\eta_k},\Gamma_{\xi_k\eta_k}) ,
\label{120}
\end{equation}
where $u_{0\xi_k\eta_k}$ and $\Gamma_{\xi_k\eta_k}$ are respectively the
projections of the point $u_0$ and of 
the loop $\Gamma$ on the plane defined by the axes $\xi_k$ and $\eta_k$,
as shown in Fig. 3.

\subsection{Factorization of polar n-complex polynomials}

A polynomial of degree $m$ of the polar n-complex variable $u$ has the form
\begin{equation}
P_m(u)=u^m+a_1 u^{m-1}+\cdots+a_{m-1} u +a_m ,
\end{equation}
where $a_l$, $l=1,...,m$, are in general polar n-complex constants.
If $a_l=\sum_{p=0}^{n-1}h_p a_{lp}$, 
and with the
notations of Eqs. (\ref{88a})-(\ref{88d}) applied for $l= 1, \cdots, m$, the
polynomial $P_m(u)$ can be written as 
\begin{eqnarray}
\lefteqn{P_m= 
e_+\left(v_+^m +\sum_{l=1}^{m}A_{l+}v_+^{m-l} \right)
+F(n)e_-\left(v_-^m +\sum_{l=1}^{m}A_{l-}v_-^{m-l} \right) \nonumber}\\
&&+\sum_{k=1}^{[(n-1)/2]}
\left[(e_k v_k+\tilde e_k\tilde v_k)^m+
\sum_{l=1}^m(e_k A_{lk}+\tilde e_k\tilde A_{lk})
(e_k v_k+\tilde e_k\tilde v_k)^{m-l} 
\right],
\label{126a}
\end{eqnarray}
where the constants $A_{l+}, A_{l-}, A_{lk}, \tilde A_{lk}$ 
are real numbers.

These relations can be written with the aid of Eqs. (\ref{e13a})-(\ref{e13b})
as 
\begin{equation}
P_m(u)=\prod_{p=1}^m (u-u_p) ,
\end{equation}
where
\begin{equation}
u_p=e_+ v_{p+}+F(n)e_-v_{p-}
+\sum_{k=1}^{[(n-1)/2]}\left(e_k v_{kp}+\tilde e_k\tilde v_{kp}\right),
\label{128d}
\end{equation}
for $ p=1,...,m$.
The roots $v_{p+}$, the roots $v_{p-}$ and, for a given $k$, the roots 
$e_k v_{kp}+\tilde e_k\tilde v_{kp}$
defined in Eq. (\ref{126a}) may be ordered arbitrarily.
This means that Eq. (\ref{128d}) gives sets of $m$ roots
$u_1,...,u_m$ of the polynomial $P_m(u)$, 
corresponding to the various ways in which the roots $v_{p+}, v_{p-}$,
$e_k v_{kp}+\tilde e_k\tilde v_{kp}$ are ordered according to $p$ in each
group. Thus, while the polar hypercomplex components in Eq. (\ref{126a}) taken
separately have unique factorizations, the polynomial $P_m(u)$ can be written
in many different ways as a product of linear factors. 

For example,
$u^2-1=(u-u_1)(u-u_2)$, where for even
$n$, $u_1=\pm e_+\pm e_-\pm  e_1\pm e_2\pm\cdots
\pm e_{n/2-1}, u_2=-u_1$, so that
there are $2^{n/2}$ independent sets of roots $u_1,u_2$
of $u^2-1$. 
It can be checked that 
$(\pm e_+\pm e_-\pm  e_1\pm e_2\pm\cdots\pm  e_{n/2-1})^2= 
e_++e_-+e_1+e_2+\cdots+e_{n/2-1}=1$.
For odd $n$, $u_1=\pm e_+\pm e_1\pm e_2\pm\cdots
\pm  e_{(n-1)/2}, u_2=-u_1$, so that there are $2^{(n-1)/2}$ independent
sets of roots $u_1,u_2$ of $u^2-1$.
It can be checked that 
$(\pm e_+\pm  e_1\pm e_2\pm\cdots\pm  e_{(n-1)/2})^2= 
e_++e_1+e_2+\cdots+e_{(n-1)/2}=1$.

\subsection{Representation of polar n-complex numbers by irreducible matrices}

The polar n-complex number 
$u$ ca be represented by the matrix
\begin{equation}
U=\left(
\begin{array}{ccccc}
x_0     &   x_1     &   x_2   & \cdots  &x_{n-1}\\
x_{n-1} &   x_0     &   x_1   & \cdots  &x_{n-2}\\
x_{n-2} &   x_{n-1} &   x_0   & \cdots  &x_{n-3}\\
\vdots  &  \vdots   &  \vdots & \cdots  &\vdots \\
x_1     &   x_2     &   x_3   & \cdots  &x_0\\
\end{array}
\right).
\label{24b}
\end{equation}
The product $u=u^\prime u^{\prime\prime}$ is
represented by the matrix multiplication $U=U^\prime U^{\prime\prime}$.
It can be shown that the irreducible form \cite{4} of the matrix $U$
in terms of matrices with real coefficients is, for even $n$,
\begin{equation}
\left(
\begin{array}{ccccc}
v_+     &     0     &     0   & \cdots  &   0   \\
0       &     v_-   &     0   & \cdots  &   0   \\
0       &     0     &     V_1 & \cdots  &   0   \\
\vdots  &  \vdots   &  \vdots & \cdots  &\vdots \\
0       &     0     &     0   & \cdots  &   V_{n/2-1}\\
\end{array}
\right)
\label{129a}
\end{equation}
and, for odd $n$, 
\begin{equation}
\left(
\begin{array}{ccccc}
v_+     &     0     &     0   & \cdots  &   0   \\
0       &     V_1   &     0   & \cdots  &   0   \\
0       &     0     &     V_2 & \cdots  &   0   \\
\vdots  &  \vdots   &  \vdots & \cdots  &\vdots \\
0       &     0     &     0   & \cdots  &   V_{(n-1)/2}\\
\end{array}
\right),
\label{129b}
\end{equation}
where 
\begin{equation}
V_k=\left(
\begin{array}{cc}
v_k           &     \tilde v_k   \\
-\tilde v_k   &     v_k          \\
\end{array}\right),
\label{130}
\end{equation}
$k=1,...,[(n-1)/2]$.
The relations between the variables $v_+, v_-, v_k, \tilde v_k$ for the
multiplication of polar n-complex numbers are
$v_+=v_+^\prime v_+^{\prime\prime},\;v_-=v_-^\prime v_-^{\prime\prime},
v_k=v_k^\prime v_k^{\prime\prime}-\tilde v_k^\prime \tilde v_k^{\prime\prime},
\;
\tilde v_k=v_k^\prime \tilde v_k^{\prime\prime}
+\tilde v_k^\prime v_k^{\prime\prime}$.

\section{Planar Hypercomplex Numbers in Even $n$ Dimensions}

\subsection{Operations with planar n-complex numbers}

A planar hypercomplex number in $n$ dimensions is determined by its $n$
components $(x_0,x_1,...,x_{n-1})$. 
The planar n-complex numbers and
their operations discussed in this paper can be represented 
by  writing the n-complex number $(x_0,x_1,...,x_{n-1})$ as  
$u=x_0+h_1x_1+h_2x_2+\cdots+h_{n-1}x_{n-1}$, where 
$h_1, h_2, \cdots, h_{n-1}$ are bases for which the multiplication rules are 
\begin{equation}
h_j h_k =(-1)^{[(j+k)/n]}h_l ,\:l=j+k-n[(j+k)/n],
\label{planar-1}
\end{equation}
for $j,k,l=0,1,..., n-1$,
where $h_0=1$. The rules for the planar bases differ from the rules for the
polar bases by the minus sign which appears when $n\leq j+k\leq 2n-2$.
The significance of the composition laws in Eq.
(\ref{planar-1}) can be understood by representing the bases $h_j, h_k$ by
points on a 
circle at the angles $\alpha_j=\pi j/n,\alpha_k=\pi k/n$, as shown in Fig. 4,
and the product $h_j h_k$ by the point of the circle at the angle 
$\pi (j+k)/n$. If $\pi\leq \pi (j+k)/n<2\pi$, the point is opposite to the
basis $h_l$ of angle $\alpha_l=\pi (j+k)/n-\pi$.

In an odd number of dimensions $n$, a transformation of coordinates according
to $x_{2l}=x^\prime_l, x_{2m-1}=-x^\prime_{(n-1)/2+m}$, 
and of the bases according to
$_{2l}=h^\prime_l, h_{2m-1}=-h^\prime_{(n-1)/2+m}, l=0,...,(n-1)/2, \;
m=1,...,(n-1)/2$,
leaves the expression of an n-complex number unchanged,
$\sum_{k=0}^{n-1}h_k x_k=\sum_{k=0}^{n-1}h^\prime_k x^\prime_k$,
and the products of the bases $h^\prime_k$ are
$h^\prime_j h^\prime_k =h^\prime_l ,\:l=j+k-n[(j+k)/n],\:\:j,k,l=0,1,..., n-1$.
Thus, the planar n-complex numbers with the rules are equivalent in an
odd number of dimensions to the polar n-complex numbers.
Therefore, in this section it will be supposed that $n$ is an even number,
unless otherwise stated.

Two n-complex numbers 
$u=x_0+h_1x_1+h_2x_2+\cdots+h_{n-1}x_{n-1}$,
$u^\prime=x^\prime_0+h_1x^\prime_1+h_2x^\prime_2+\cdots+h_{n-1}x^\prime_{n-1}$ 
are equal if and only if $x_j=x^\prime_j, j=0,1,...,n-1$.
The sum of the n-complex numbers $u$
and
$u^\prime$ 
is
\begin{equation}
u+u^\prime=x_0+x^\prime_0+h_1(x_1+x^\prime_1)+\cdots
+h_{n-1}(x_{n-1} +x^\prime_{n-1}) .
\label{planar-2}
\end{equation}
The product of the numbers $u, u^\prime$ is 
\begin{equation}
\begin{array}{l}
uu^\prime=x_0 x_0^\prime -x_1x_{n-1}^\prime-x_2
x_{n-2}^\prime-x_3x_{n-3}^\prime 
-\cdots-x_{n-1}x_1^\prime\\
+h_1(x_0 x_1^\prime+x_1x_0^\prime-x_2x_{n-1}^\prime-x_3x_{n-2}^\prime
-\cdots-x_{n-1} x_2^\prime) \\
+h_2(x_0 x_2^\prime+x_1x_1^\prime+x_2x_0^\prime-x_3x_{n-1}^\prime
-\cdots-x_{n-1} x_3^\prime) \\
\vdots\\
+h_{n-1}(x_0
x_{n-1}^\prime+x_1x_{n-2}^\prime+x_2x_{n-3}^\prime+x_3x_{n-4}^\prime 
+\cdots+x_{n-1} x_0^\prime).
\end{array}
\label{planar-3}
\end{equation}
The product $uu^\prime$ can be written as
\begin{equation}
uu^\prime=\sum_{k=0}^{n-1}h_k\sum_{l=0}^{n-1}(-1)^{[(n-k-1+l)/n]}x_l
x^\prime_{k-l+n[(n-k-1+l)/n]}. 
\label{planar-3a}
\end{equation}
If $u,u^\prime,u^{\prime\prime}$ are n-complex numbers, the multiplication 
is associative,
$(uu^\prime)u^{\prime\prime}=u(u^\prime u^{\prime\prime})$,
and commutative,
$u u^\prime=u^\prime u$ ,
because the product of the bases, defined in Eq. (\ref{planar-1}), is
associative and 
commutative. 

The inverse of the planar n-complex number $u$ is the n-complex number
$u^\prime$ having the property that $uu^\prime=1$.
This equation has a solution provided that the corresponding determinant
$\nu$ is not equal to zero, $\nu\not=0$.
For planar n-complex numbers $\nu\geq 0$, and the quantity 
$\rho=\nu^{1/n}$
will be called amplitude of the
n-complex number $u$.
It can be shown that 
\begin{equation}
\nu=\prod_{k=1}^{n/2}\rho_k^2,
\end{equation}
where
\begin{equation}
\rho_k^2=v_k^2+\tilde v_k^2,
\end{equation}
\begin{equation}
v_k=\sum_{p=0}^{n-1}x_p\cos\left(\frac{\pi (2k-1)p}{n}\right),
\end{equation}
\begin{equation}
\tilde v_k=\sum_{p=0}^{n-1}x_p\sin\left(\frac{\pi (2k-1)p}{n}\right).
\end{equation}
Thus, a planar n-complex number has an inverse
unless it lies on one of the nodal hypersurfaces 
$\rho_1=0$, or $\rho_2=0$, or ... or $\rho_{n/2}=0$.
It can also be shown that
\begin{equation}
d^2=\frac{2}{n}\sum_{k=1}^{n/2}\rho_k^2.
\end{equation}
From this relation it results that if the product of two n-complex numbers is
zero, $uu^\prime=0$, then 
$\rho_k\rho_k^\prime=0, k=1,...,n/2$, which means that either $u=0$, or
$u^\prime=0$, or $u, u^\prime$ belong to orthogonal hypersurfaces in such a way
that the afore-mentioned products of components should be equal to zero.

\subsection{Geometric representation of planar n-complex numbers}

The planar n-complex number $x_0+h_1x_1+h_2x_2+\cdots+h_{n-1}x_{n-1}$
can be represented by 
the point $A$ of coordinates $(x_0,x_1,...,x_{n-1})$. 
If $O$ is the origin of the n-dimensional space,  the
distance from the origin $O$ to the point $A$ of coordinates
$(x_0,x_1,...,x_{n-1})$ has the expression written in Eq. (\ref{10}).
The quantity $d$ will be called now modulus of the planar n-complex number 
$u=x_0+h_1x_1+h_2x_2+\cdots+h_{n-1}x_{n-1}$. The modulus of an n-complex number
$u$ will be designated by $d=|u|$. The quantity $\rho=\nu^{1/n}$
will be called amplitude of the n-complex number $u$. 

The exponential and trigonometric forms of the n-complex number $u$ can be
obtained conveniently in a rotated system of axes defined by the transformation
\begin{equation}
v_k= \sqrt{n/2}\xi_k , \tilde v_k= \sqrt{n/2}\eta_k, 
\end{equation}
for $k=1,...,n/2$.
This transformation from the coordinates $x_0,...,x_{n-1}$ to the variables
$\xi_k, \eta_k$ is unitary.

The position of the point $A$ of coordinates $(x_0,x_1,...,x_{n-1})$ can be
also described with the aid of the distance $d$, Eq. (\ref{10}), and of $n-1$
angles defined further. Thus, in the plane of the axes $v_k,\tilde v_k$, the
radius $\rho_k$ and the azimuthal angle 
$\phi_k$ can be introduced by the relations 
\begin{equation}
\cos\phi_k=v_k/\rho_k,\:\sin\phi_k=\tilde v_k/\rho_k, 
\end{equation}
for
$0\leq \phi_k<2\pi, k=1,...,n/2$, 
so that there are $n/2$ azimuthal angles.
If the projection of the point $A$ on the plane of the axes $v_k,\tilde v_k$ is
$A_k$, 
and the projection of the point $A$ on the 4-dimensional space defined by the
axes $v_1, \tilde v_1, v_k,\tilde v_k$ is $A_{1k}$, the angle $\psi_{k-1}$
between the line 
$OA_{1k}$ and the 2-dimensional plane defined by the axes $v_k,\tilde v_k$ is 
\begin{equation}
\tan\psi_{k-1}=\rho_1/\rho_k, 
\end{equation}
where $0\leq\psi_k\leq\pi/2, k=2,...,n/2$,
so that there are $n/2-1$ planar angles.
Thus, the position of
the point $A$ is described by the distance
$d$, by $n/2$ azimuthal angles and by $n/2-1$ planar angles. The variables
$\nu, \rho, 
\rho_k, \tan\psi_k$ are multiplicative and the azimuthal angles $\phi_k$ are
additive upon the multiplication of polar n-complex numbers.

\subsection{The planar n-dimensional cosexponential functions}

The exponential function of the planar n-complex variable $u$ can be defined by
the series
$\exp u = 1+u+u^2/2!+u^3/3!+\cdots$ . 
It can be checked by direct multiplication of the series that
$\exp(u+u^\prime)=\exp u \cdot \exp u^\prime$ ,
so that  
$\exp u=\exp x_0 \cdot \exp (h_1x_1) \cdots \exp (h_{n-1}x_{n-1})$.

It can be seen with the aid of the representation in Fig. 4 that 
\begin{equation}
h_k^{n+p}=-h_k^p, 
\end{equation}
for $p$ integer, $k=1,...,n-1$.
For $k$ even, $e^{h_k y}$ can be written as
\begin{equation}
e^{h_k y}=\sum_{p=0}^{n-1}(-1)^{[kp/n]}h_{kp-n[kp/n]}g_{np}(y),
\label{planar-28bx}
\end{equation}
where $g_{np}$ are the polar n-dimensional cosexponential
functions.  For odd $k$,  $e^{h_k y}$ is
\begin{equation}
e^{h_k y}=\sum_{p=0}^{n-1}(-1)^{[kp/n]}h_{kp-n[kp/n]}f_{np}(y),
\label{planar-28b}
\end{equation}
where the functions $f_{nk}$, which will be
called planar cosexponential functions in $n$ dimensions, are
\begin{equation}
f_{nk}(y)=\sum_{p=0}^\infty (-1)^p \frac{y^{k+pn}}{(k+pn)!},
\label{planar-29}
\end{equation}
for $ k=0,1,...,n-1$.

The planar cosexponential functions  of even index $k$ are
even functions, $f_{n,2l}(-y)=f_{n,2l}(y)$, 
and the planar cosexponential functions of odd index 
are odd functions, $f_{n,2l+1}(-y)=-f_{n,2l+1}(y)$, $l=0,...,n/2-1$ . 

The planar n-dimensional cosexponential function $f_{nk}(y)$ is related to the
polar n-dimensional cosexponential function $g_{nk}(y)$ by the relation 
\begin{equation}
f_{nk}(y)=e^{-i\pi k/n}g_{nk}\left(e^{i\pi/n}y\right), 
\label{planar-30a}
\end{equation}
for $k=0,...,n-1$.
The expression of the planar n-dimensional cosexponential functions is then
\begin{equation}
f_{nk}(y)=\frac{1}{n}\sum_{l=1}^{n}\exp\left[y\cos\left(\frac{\pi
(2l-1)}{n}\right) \right]
\cos\left[y\sin\left(\frac{\pi (2l-1)}{n}\right)-\frac{\pi (2l-1)k}{n}\right],
\label{planar-30}
\end{equation}
$k=0,1,...,n-1$. 
The planar cosexponential function defined in Eq. (\ref{planar-29}) has the
expression 
given in Eq. (\ref{planar-30}) for any natural value of $n$, this result not
being 
restricted to even values of $n$.

It can be shown from Eq. (\ref{planar-30}) that
\begin{equation}
\sum_{k=0}^{n-1}f_{nk}^2(y)=\frac{1}{n}\sum_{l=1}^{n}\exp\left[2y\cos\left(
\frac{\pi (2l-1)}{n}\right)\right].
\label{planar-34a}
\end{equation}
It can be seen that the right-hand side of Eq. (\ref{planar-34a}) does not
contain 
oscillatory terms. If $n$ is a multiple of 4, it can be shown by replacing $y$
by $iy$ in Eq. (\ref{planar-34a}) that
\begin{equation}
\sum_{k=0}^{n-1}(-1)^kf_{nk}^2(y)=\frac{4}{n}
\sum_{l=1}^{n/4}\cos\left[2y\cos\left(\frac{\pi (2l-1)}{n}\right)\right],
\label{planar-34b}
\end{equation}
which does not contain exponential terms.

For odd $n$, the planar n-dimensional cosexponential function $f_{nk}(y)$ is
related to the n-dimensional cosexponential function $g_{nk}(y)$ also by the
relation  
\begin{equation}
f_{nk}(y)=(-1)^kg_{nk}(-y),
\end{equation}
as can be seen by comparing the series for the two classes of functions.

Addition theorems for the planar n-dimensional cosexponential functions can be
obtained from the relation $\exp h_1(y+z)=\exp h_1 y \cdot\exp h_1 z $, by
substituting the expression of the exponentials as given 
by $e^{h_1 y}=\sum_{p=0}^{n-1} h_p f_{np}(y)$,
\begin{eqnarray}
\lefteqn{f_{nk}(y+z)=f_{n0}(y)f_{nk}(z)+f_{n1}(y)f_{n,k-1}(z)+\cdots
+f_{nk}(y)f_{n0}(z)\nonumber}\\
&&-f_{n,k+1}(y)f_{n,n-1}(z)-f_{n,k+2}(y)f_{n,n-2}(z)
-\cdots-f_{n,n-1}(y)f_{n,k+1}(z) ,
\label{planar-35a}
\end{eqnarray}
for $k=0,1,...,n-1$.
It can also be shown that
\begin{equation}
\left\{f_{n0}(y)+h_1f_{n1}(y)+\cdots+h_{n-1}f_{n,n-1}(y)\right\}^l
=f_{n0}(ly)+h_1f_{n1}(ly)+\cdots+h_{n-1}f_{n,n-1}(ly).
\label{planar-37b}
\end{equation}

The planar n-dimensional cosexponential functions are solutions of the
$n^{\rm th}$-order differential equation
\begin{equation}
\frac{{\rm d}^n\zeta}{{\rm d}u^n}=-\zeta ,
\label{planar-44}
\end{equation}
whose solutions are of the form
$\zeta(u)=A_0f_{n0}(u)+A_1f_{n1}(u)+\cdots+A_{n-1}f_{n,n-1}(u)$.
It can be checked that the derivatives of the planar cosexponential functions
are related by
\begin{equation}
\frac{df_{n0}}{du}=-f_{n,n-1}, \:
\frac{df_{n1}}{du}=f_{n0}, \:...,
\frac{df_{n,n-2}}{du}=f_{n,n-3} ,
\frac{df_{n,n-1}}{du}=f_{n,n-2} .
\label{planar-45}
\end{equation}

For $n=2$, the planar cosexponential functions are $f_{20}(y)=\cos y$ and
$f_{21}(y)=\sin y$.

\subsection{Exponential and trigonometric forms of planar n-complex numbers}

In order to obtain the exponential and trigonometric forms of planar n-complex
numbers, a new set of hypercomplex bases will be introduced 
by the relations 
\begin{equation}
e_k=\frac{2}{n}\sum_{p=0}^{n-1}h_p\cos\left( \frac{\pi(2k-1)p}{n}\right),
\label{planar-e13a}
\end{equation}
\begin{equation}
\tilde e_k=\frac{2}{n}\sum_{p=0}^{n-1}h_p
\sin\left( \frac{\pi(2k-1)p}{n}\right),
\label{planar-e13b}
\end{equation}
for $k=1,...,n/2$.
The multiplication relations for the new hypercomplex bases are
\begin{eqnarray}
e_k^2=e_k, \tilde e_k^2=-e_k, e_k \tilde e_k=\tilde e_k , e_ke_l=0, e_k\tilde
e_l=0, \tilde e_k\tilde e_l=0, k\not=l, 
\label{planar-e12a}
\end{eqnarray}
for $k,l=1,...,n/2$.
It can be shown that
\begin{equation}
x_0+h_1x_1+\cdots+h_{n-1}x_{n-1}= 
\sum_{k=1}^{n/2}(e_k v_k+\tilde e_k \tilde v_k).
\end{equation}

The exponential form of the planar n-complex number $u$ is
\begin{eqnarray}
u=\rho\exp\left\{\sum_{p=1}^{n-1}h_p\left[
-\frac{2}{n}\sum_{k=2}^{n/2}
\cos\left(\frac{\pi (2k-1)p}{n}\right)\ln\tan\psi_{k-1}
\right]
+\sum_{k=1}^{n/2}\tilde e_k\phi_k
\right\},
\label{planar-50a}
\end{eqnarray}
where the amplitude is
\begin{equation}
\rho=\left(\rho_1^2\cdots\rho_{n/2}^2\right)^{1/n}.
\end{equation}

The trigonometric form of the planar n-complex number $u$ is
\begin{eqnarray}
\lefteqn{u=d
\left(\frac{n}{2}\right)^{1/2}
\left(1+\frac{1}{\tan^2\psi_1}+\frac{1}{\tan^2\psi_2}+\cdots
+\frac{1}{\tan^2\psi_{n/2-1}}\right)^{-1/2}\nonumber}\\
&&\left(e_1+\sum_{k=2}^{n/2}\frac{e_k}{\tan\psi_{k-1}}\right)
\exp\left(\sum_{k=1}^{n/2}\tilde e_k\phi_k\right).
\label{planar-52a}
\end{eqnarray}

\subsection{Elementary functions of a planar n-complex variable}

The logarithm $u_1$ of the planar n-complex number $u$, $u_1=\ln u$, can be
defined as the solution of the equation
$u=e^{u_1}$ .
The logarithm exists as a planar n-complex function with real components for
all values of $x_0,...,x_{n-1}$ for which $\rho\not=0$. 
The expression of the logarithm is
\begin{eqnarray}
\ln u=\ln \rho+\sum_{p=1}^{n-1}h_p\left[
-\frac{2}{n}\sum_{k=2}^{n/2}
\cos\left(\frac{\pi (2k-1)p}{n}\right)\ln\tan\psi_{k-1}
\right]+\sum_{k=1}^{n/2}\tilde e_k\phi_k.
\label{planar-56a}
\end{eqnarray}
The function $\ln u$ is multivalued because of the presence of the terms 
$\tilde e_k\phi_k$.

The power function $u^m$ of the planar n-complex variable $u$ can be defined
for real values of $m$ as $u^m=e^{m\ln u}$ .
It can be shown that
\begin{equation}
u^m=\sum_{k=1}^{n/2}
\rho_k^m(e_k\cos m\phi_k+\tilde e_k\sin m\phi_k).
\label{planar-59a}
\end{equation}
The power function is multivalued unless $m$ is an integer.

\subsection{Power series of planar n-complex numbers}

A power series of the planar n-complex variable $u$ is a series of the form
\begin{equation}
a_0+a_1 u + a_2 u^2+\cdots +a_l u^l+\cdots .
\label{planar-83}
\end{equation}
Using the inequality 
\begin{equation}
|u^\prime u^{\prime\prime}|\leq
\sqrt{n/2}|u^\prime||u^{\prime\prime}| ,
\label{planar-83b}
\end{equation}
which replaces the relation of equality
extant for 2-dimensional complex numbers, it can be shown that
the series (\ref{planar-83}) is absolutely convergent for 
$|u|<c$,
where 
$c=\lim_{l\rightarrow\infty} |a_l|/\sqrt{n/2}|a_{l+1}|$ .

The convergence of the series (\ref{planar-83}) can be also studied with the
aid of 
the formulas (\ref{planar-59a}), which for integer values of $m$ are
valid for any values of $x_0,...,x_{n-1}$.
If 
$a_l=\sum_{p=0}^{n-1}h_p a_{lp}$, and
\begin{equation}
A_{lk}=\sum_{p=0}^{n-1} a_{lp}\cos\frac{\pi (2k-1)p}{n},
\label{planar-88a}
\end{equation}
\begin{equation}
\tilde A_{lk}=\sum_{p=0}^{n-1} a_{lp}\sin\frac{\pi (2k-1)p}{n},
\label{planar-88b}
\end{equation}
for $k=1,...,n/2$,
the series (\ref{planar-83}) can be written as
\begin{equation}
\sum_{l=0}^\infty \sum_{k=1}^{n/2}
(e_k A_{lk}+\tilde e_k\tilde A_{lk})(e_k v_k+\tilde e_k\tilde v_k)^l.
\label{planar-89a}
\end{equation}

The series in Eq. (\ref{planar-83}) is absolutely convergent for 
\begin{equation}
\rho_k<c_k, 
\end{equation}
for $k=1,..., n/2$, where
\begin{equation}
c_k=\lim_{l\rightarrow\infty} \frac
{\left[A_{lk}^2+\tilde A_{lk}^2\right]^{1/2}}
{\left[A_{l+1,k}^2+\tilde A_{l+1,k}^2\right]^{1/2}} .
\end{equation}
These relations show that the region of convergence of the series
(\ref{planar-83}) is an n-dimensional cylinder.

\subsection{Analytic functions of planar n-complex variables}

If the n-complex function $f(u)$
of the planar n-complex variable $u$ is written in terms of 
the real functions $P_k(x_0,...,x_{n-1}), k=0,1,...,n-1$ of the real
variables $x_0,x_1,...,x_{n-1}$ as 
\begin{equation}
f(u)=\sum_{k=0}^{n-1}h_kP_k(x_0,...,x_{n-1}),
\label{planar-h93}
\end{equation}
then relations of equality 
exist between the partial derivatives of the functions $P_k$,
\begin{equation}
\frac{\partial P_k}{\partial x_0} = \frac{\partial P_{k+1}}{\partial x_1} 
=\cdots=\frac{\partial P_{n-1}}{\partial x_{n-k-1}} 
= -\frac{\partial P_0}{\partial x_{n-k}}=\cdots
=-\frac{\partial P_{k-1}}{\partial x_{n-1}}, 
\label{planar-h95}
\end{equation}
for $k=0,1,...,n-1$.
The relations (\ref{planar-h95}) are analogous to the Riemann relations
for the real and imaginary components of a complex function. 
It can be shown from Eqs. (\ref{planar-h95}) that the components $P_k$ fulfil
the 
second-order equations
\begin{eqnarray}
\lefteqn{\frac{\partial^2 P_k}{\partial x_0\partial x_l}
=\frac{\partial^2 P_k}{\partial x_1\partial x_{l-1}}
=\cdots=
\frac{\partial^2 P_k}{\partial x_{[l/2]}\partial x_{l-[l/2]}}}\nonumber\\
&&=-\frac{\partial^2 P_k}{\partial x_{l+1}\partial x_{n-1}}
=-\frac{\partial^2 P_k}{\partial x_{l+2}\partial x_{n-2}}
=\cdots
=-\frac{\partial^2 P_k}{\partial x_{l+1+[(n-l-2)/2]}
\partial x_{n-1-[(n-l-2)/2]}} ,\nonumber\\
&&
\label{planar-96}
\end{eqnarray}
for $k,l=0,1,...,n-1$.

\subsection{Integrals of planar n-complex functions}

The singularities of planar n-complex functions arise from terms of the form
$1/(u-u_0)^m$, with $m>0$. Functions containing such terms are singular not
only at $u=u_0$, but also at all points of the hypersurfaces
passing through $u_0$ and which are parallel to the nodal hypersurfaces. 

The integral of a planar n-complex function between two points $A, B$ along a
path 
situated in a region free of singularities is independent of path, which means
that the integral of an analytic function along a loop situated in a region
free of singularities is zero,
\begin{equation}
\oint_\Gamma f(u) du = 0,
\label{planar-111}
\end{equation}
where it is supposed that a surface $\Sigma$ spanning 
the closed loop $\Gamma$ is not intersected by any of
the hypersurfaces associated with the
singularities of the function $f(u)$. Using the expression, 
Eq. (\ref{planar-h93}), for $f(u)$ and the fact that 
$du=\sum_{k=0}^{n-1}h_k dx_k$, 
the explicit form of the integral in Eq. (\ref{planar-111}) is
\begin{equation}
\oint _\Gamma f(u) du = \oint_\Gamma
\sum_{k=0}^{n-1}h_k\sum_{l=0}^{n-1}(-1)^{[(n-k-1+l)/n]}
P_l dx_{k-l+n[(n-k-1+l)/n]}.
\end{equation}
If the functions $P_k$ are regular on a surface $\Sigma$
spanning the loop $\Gamma$,
the integral along the loop $\Gamma$ can be transformed in an integral over the
surface $\Sigma$ of terms of the form
$\partial P_l/\partial x_{k-m+n[(n-k+m-1)/n]}$
$-(-1)^s \partial P_m/\partial x_{k-l+n[(n-k+l-1)/n]}$, where 
$s=[(n-k+m-1)/n]-[(n-k+l-1)/n]$.
The integrals of these terms are equal to zero by Eqs. (\ref{planar-h95}), and
this 
proves Eq. (\ref{planar-111}). 

The quantity $du/(u-u_0)$ is
\begin{eqnarray}
\frac{du}{u-u_0}=
\frac{d\rho}{\rho}+\sum_{p=1}^{n-1}h_p\left[
-\frac{2}{n}\sum_{k=2}^{n/2}
\cos\left(\frac{2\pi kp}{n}\right)d\ln\tan\psi_{k-1}\right]
+\sum_{k=1}^{n/2}\tilde e_kd\phi_k.
\label{planar-114a}
\end{eqnarray}
Since $\rho$ and $\ln(\tan\psi_{k-1})$ are singlevalued variables, it follows
that 
$\oint_\Gamma d\rho/\rho =0$, and 
$\oint_\Gamma d(\ln\tan\psi_{k-1})=0$.
On the other hand, since $\phi_k$ are cyclic variables, they may give
contributions to the integral around the closed loop $\Gamma$.

If $f(u)$ is an analytic function of a polar n-complex variable which can be
expanded in a 
series which holds on the curve
$\Gamma$ and on a surface spanning $\Gamma$, then 
\begin{equation}
\oint_\Gamma \frac{f(u)du}{u-u_0}=
2\pi f(u_0)\sum_{k=1}^{n/2}\tilde e_k 
\;{\rm int}(u_{0\xi_k\eta_k},\Gamma_{\xi_k\eta_k}) .
\label{planar-120}
\end{equation}

\subsection{Factorization of planar n-complex polynomials}

A polynomial of degree $m$ of the planar n-complex variable $u$ has the form
\begin{equation}
P_m(u)=u^m+a_1 u^{m-1}+\cdots+a_{m-1} u +a_m ,
\end{equation}
where $a_l$, $l=1,...,m$, are in general planar n-complex constants.
If $a_l=\sum_{p=0}^{n-1}h_p a_{lp}$, and with the
notations of Eqs. (\ref{planar-88a})-(\ref{planar-88b}) applied for $l= 1,
\cdots, m$, the 
polynomial $P_m(u)$ can be written as 
\begin{eqnarray}
P_m= \sum_{k=1}^{n/2}\left\{(e_k v_k+\tilde e_k\tilde v_k)^m+
\sum_{l=1}^m(e_k A_{lk}+\tilde e_k\tilde A_{lk})
(e_k v_k+\tilde e_k\tilde v_k)^{m-l} 
\right\},
\label{planar-126a}
\end{eqnarray}
where the constants $A_{lk}, \tilde A_{lk}$ are real numbers.

These relations can be written with the aid of Eqs. (\ref{planar-e13a}) and
(\ref{planar-e13b}) as
\begin{equation}
P_m(u)=\prod_{p=1}^m (u-u_p) ,
\end{equation}
where
\begin{equation}
u_p=\sum_{k=1}^{n/2}\left(e_k v_{kp}
+\tilde e_k\tilde v_{kp}\right), 
\label{planar-129b}
\end{equation}
for $p=1,...,m$.
For a given $k$, the roots  
$e_k v_{kp}+\tilde e_k\tilde v_{kp}$
defined in Eq. (\ref{planar-126a}) may be ordered arbitrarily. This means that
Eq. 
(\ref{planar-129b}) gives sets of $m$ roots 
$u_1,...,u_m$ of the polynomial $P_m(u)$, 
corresponding to the various ways in which the roots $e_k v_{kp}+\tilde
e_k\tilde v_{kp}$ are ordered according to $p$ for each value of $k$. 
Thus, while the planar hypercomplex components in Eq. (\ref{planar-126a}) taken
separately have 
unique factorizations, the polynomial $P_n(u)$ can be written in many different
ways as a product of linear factors. 

For example, $u^2+1=(u-u_1)(u-u_2)$, where 
$u_1=\pm \tilde e_1\pm\tilde e_2\pm\cdots\pm \tilde e_{n/2}, u_2=-u_1$, so that
there are $2^{n/2-1}$ independent sets of roots $u_1,u_2$
of $u^2+1$. It can be checked
that $(\pm \tilde e_1\pm\tilde e_2\pm\cdots\pm \tilde e_{n/2})^2=
-e_1-e_2-\cdots-e_{n/2}=-1$.

\subsection{Representation of planar n-complex numbers by irreducible matrices}

The planar n-complex number $u$ ca be represented by the matrix
\begin{equation}
U=\left(
\begin{array}{ccccc}
x_0      &   x_1      &   x_2    & \cdots  &x_{n-1}\\
-x_{n-1} &   x_0      &   x_1    & \cdots  &x_{n-2}\\
-x_{n-2} &  - x_{n-1} &   x_0    & \cdots  &x_{n-3}\\
\vdots   &  \vdots    &  \vdots  & \cdots  &\vdots \\
-x_1     &   -x_2     &  - x_3   & \cdots  &x_0\\
\end{array}
\right).
\label{planar-24a}
\end{equation}
The product $u=u^\prime u^{\prime\prime}$ is, for even $n$,
represented by the matrix multiplication $U=U^\prime U^{\prime\prime}$.
It can be shown that the irreducible form \cite{4} of the matrix $U$,
in terms of matrices with real coefficients, is
\begin{equation}
\left(
\begin{array}{ccccc}
v_+     &         0   & \cdots  &   0   \\
0       &         V_1 & \cdots  &   0   \\
\vdots  &      \vdots & \cdots  &\vdots \\
0       &         0   & \cdots  &   V_{n/2}\\
\end{array}
\right),
\label{planar-129a}
\end{equation}
where 
\begin{equation}
V_k=\left(
\begin{array}{cc}
v_k           &     \tilde v_k   \\
-\tilde v_k   &     v_k          \\
\end{array}\right), 
\label{planar-130}
\end{equation}
for $ k=1,...,n/2$.
The relations between the variables $v_k, \tilde v_k$ for the
multiplication of planar n-complex numbers are
$v_k=v_k^\prime v_k^{\prime\prime}-\tilde v_k^\prime \tilde
v_k^{\prime\prime},\; 
\tilde v_k=v_k^\prime \tilde v_k^{\prime\prime}+\tilde v_k^\prime
v_k^{\prime\prime}$.

\section{Conclusions}

The polar and planar n-complex numbers described in this paper have a
geometric representation based on modulus, amplitude and angular variables.
The n-complex numbers have exponential and trigonometric forms, which can be
expressed with the aid of geometric variables. The exponential function of
an n-complex variable can be developed in terms of the cosexponential
functions.  The n-complex functions defined by series of powers are analytic,
and the partial derivatives of the real components of n-complex functions are
closely related. The integrals of n-complex functions are independent of path
in regions where the functions are regular. The fact that the exponential form
of the n-complex numbers depends on the cyclic azimuthal variables leads to the
concept of pole and residue for n-complex integrals on closed paths. The
polynomials of polar n-complex variables can be written as products of linear
or quadratic factors, and the polynomials of planar n-complex variables can be
written as products of linear factors.

\newpage

FIGURE CAPTIONS\\

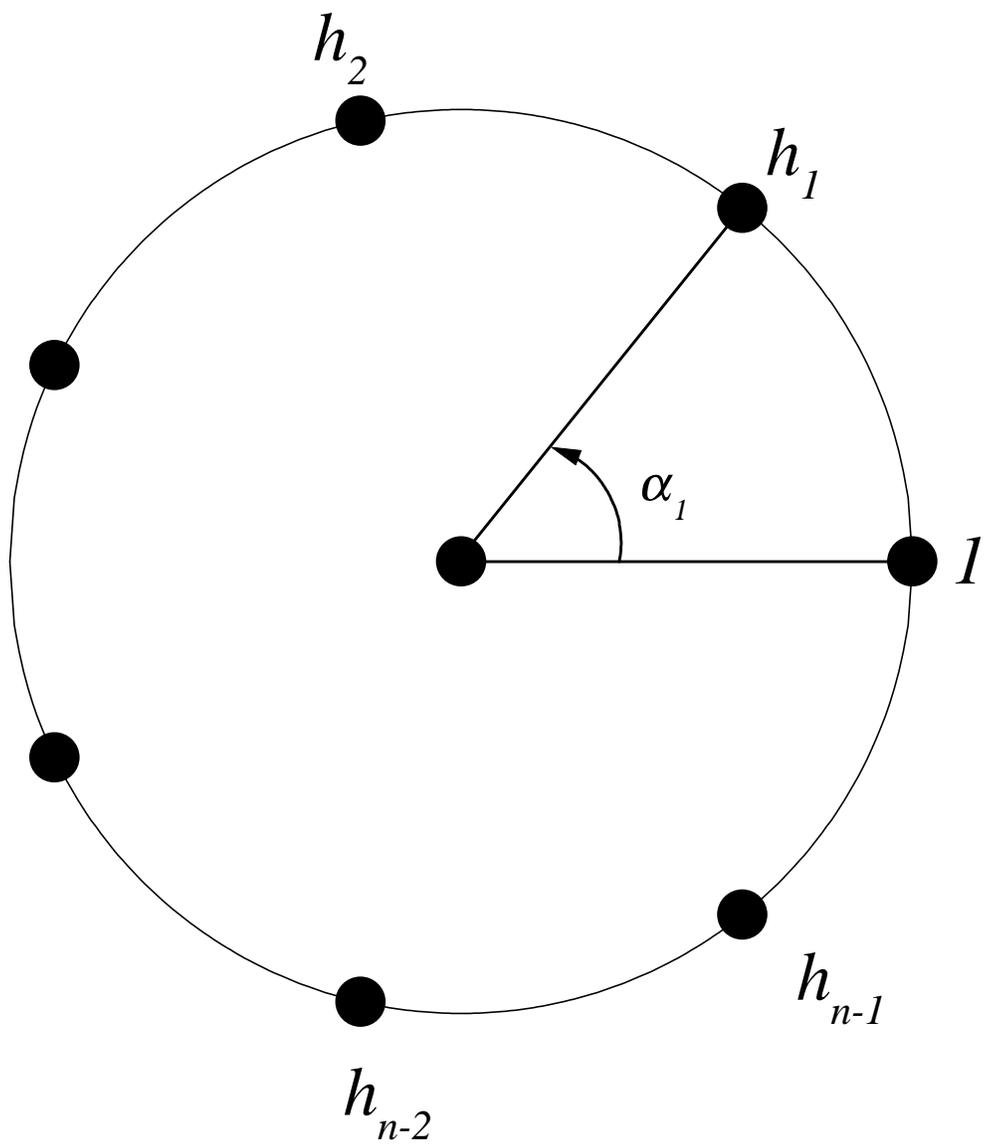
Fig. 1. Representation of the polar n-complex bases $1, h_1,...,h_{n-1}$
by points on a circle at the angles $\alpha_k=2\pi k/n$.\\

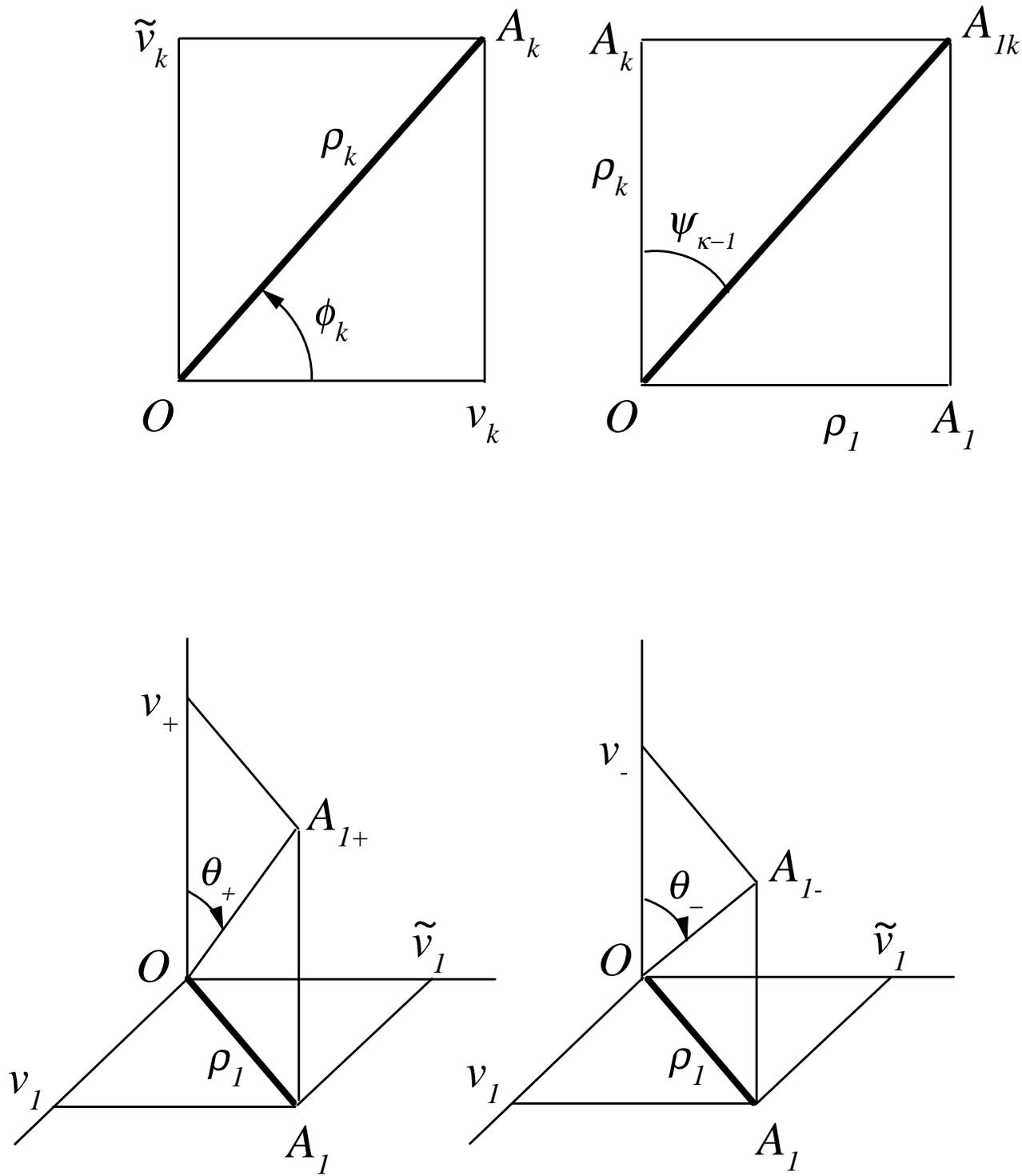
Fig. 2. Angular variables for the description of n-complex numbers.\\

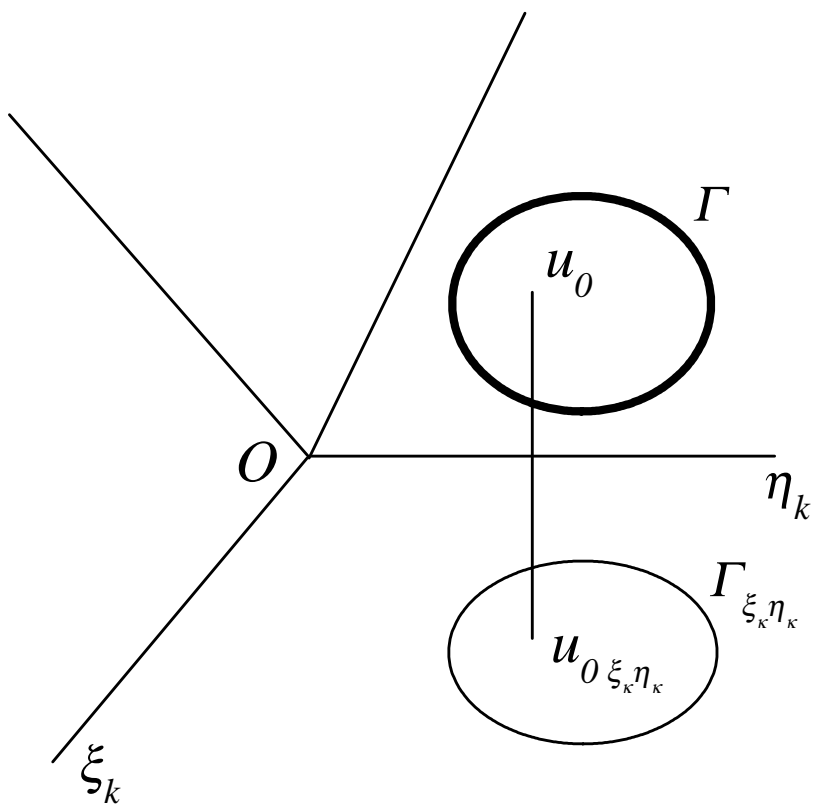
Fig. 3. Integration path $\Gamma$ and pole $u_0$, and their projections
$\Gamma_{\xi_k\eta_k}$ and $u_{0\xi_k\eta_k}$ on the plane $\xi_k \eta_k$.\\ 

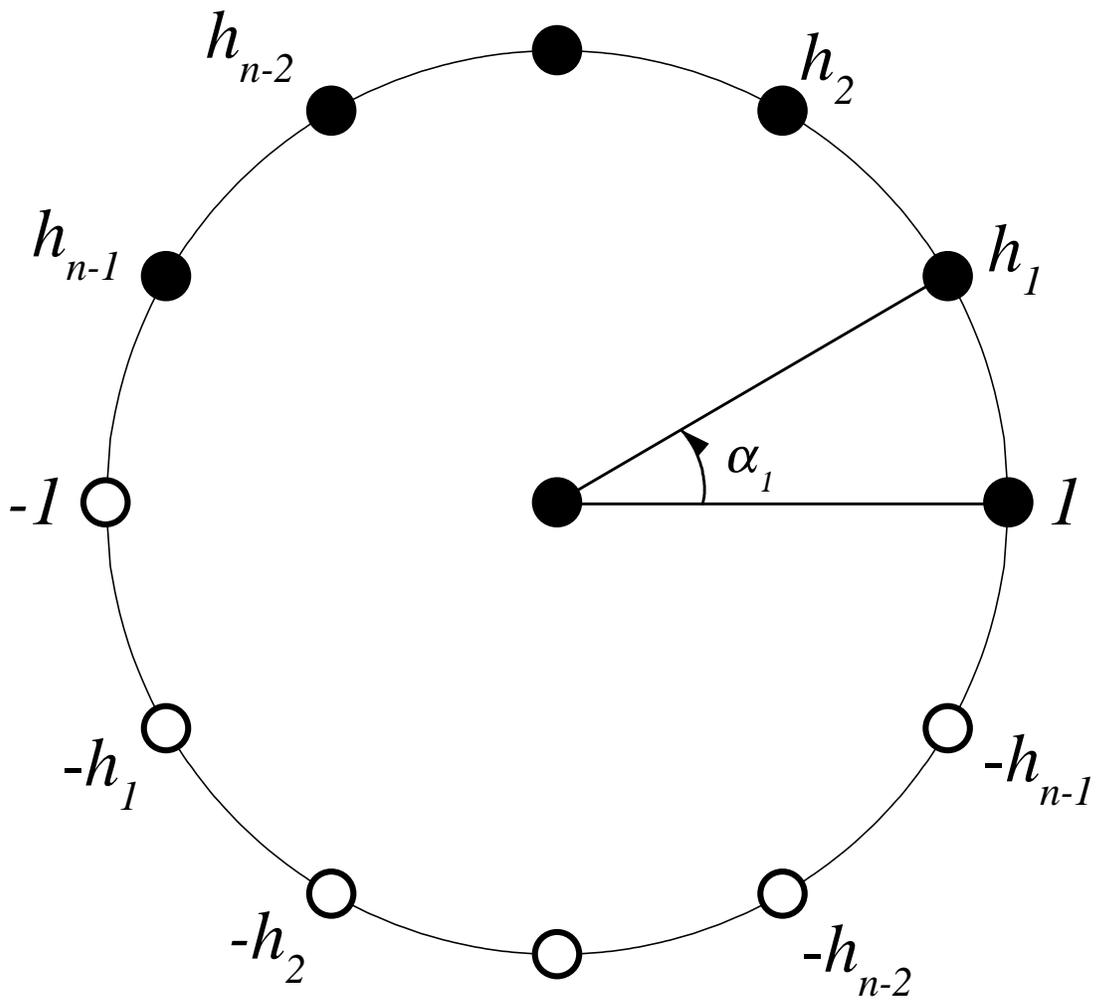
Fig. 4. Representation of the planar n-complex bases $1, h_1,...,h_{n-1}$
by points on a circle at the angles $\alpha_k=\pi k/n$.

\end{document}